\documentclass[a4paper,12pt]{article}
\usepackage{amssymb}
\usepackage{amsmath}
\usepackage{amsthm}
\usepackage{mathrsfs}
\def\dj{d\kern-.30em\raise1.25ex\vbox{\hrule width .3em height .03em}}
\def\Dj{\rlap{\kern-.70em\raise0.75ex\vbox{\hrule width .3em height .03em}}} 

\newtheorem{theorem}{Theorem}[section]
\newtheorem{proposition}[theorem]{Proposition}
\newtheorem{lemma}[theorem]{Lemma}

\newcommand{\bfrak }{\mathfrak{b}}

\newcommand{\cF}{\mathcal{F}}

\newcommand{\C }{\mathbb{C}}

\newcommand{\dif }{\mathrm{d}}

\newcommand{\gfrak}{\mathfrak{g}}

\newcommand{\gr}{\mathrm{gr}}
\newcommand{\hfrak}{\mathfrak{h}}

\newcommand{\Hom}{\mathrm{Hom}}

\newcommand{\id}{\mathrm{Id}}

\newcommand{\im}{\mathrm{Im}}

\newcommand{\ke }{\mathrm{Ker} }

\newcommand{\lfrak}{\mathfrak{l}}

\newcommand{\N }{\mathbb{N}}
\newcommand{\nfrak}{\mathfrak{n}}

\newcommand{\ot }{\otimes }

\newcommand{\pfrak}{\mathfrak{p}}

\newcommand{\U }{U}
\newcommand{\ufrak }{\mathfrak{u}}
\newcommand{\ug }{U(\mathfrak{g})}

\newcommand{\uqg }{U_q(\mathfrak{g})}

\newcommand{\unsm }{U(\mathfrak{n}_S^-)}
\newcommand{\uusm }{U(\mathfrak{u}_S^-)}

\newcommand{\wlat }{P}

\newcommand{\wurz }{\varDelta }
\newcommand{\vep }{\varepsilon }

\newcommand{\Z}{\mathbb{Z}}

\begin{document}
\title{On the Bernstein-Gelfand-Gelfand resolution for Kac-Moody algebras and quantized  
       enveloping algebras}

\author{Istv\'an Heckenberger 
        and Stefan Kolb\footnote{Supported by the German Research Foundation (DFG)}}

\date{May 16, 2006}

\maketitle

\begin{abstract}
 A Bernstein-Gelfand-Gelfand resolution for arbitrary Kac-Moody algebras and arbitrary
 subsets of the set of simple roots is proven. Moreover, quantum group analogs of the
 Bernstein-Gelfand-Gelfand resolution for symmetrizable Kac-Moody algebras are established.
 For quantized enveloping algebras with fixed deformation parameter $q\in\C\setminus\{0\}$
 exactness is proven for all $q$ which are not a root of unity.
\end{abstract}

\section{Introduction}

In this note we consider the Bernstein-Gelfand-Gelfand resolution (BGG-resolution)
of a maximal integrable highest weight module over an arbitrary Kac-Moody algebra 
$\gfrak$ and its generalization for a standard parabolic subalgebra $\pfrak_S\subset\gfrak$
(which we will abbreviate BGGL-resolution). These resolutions, which for symmetrizable 
$\gfrak$ are classical
results in representation theory by now, lead to conceptual insight into the Weyl-Kac 
character
formula and the $\nfrak^-$-homology results of Kostant and Garland-Lepowsky.

In the finite case the BGG-resolution was introduced in \cite{a-BGG75} and extended
to parabolic subalgebras by J. Lepowsky in \cite{a-Lep77}. Generalizations for 
symmetrizable Kac-Moody algebras were established soon after \cite{a-GarLep76}, \cite{a-RCW82}.
Finally S.~Kumar proved exactness of the BGG-complex for arbitrary Kac-Moody algebras
\cite{a-Kumar90} using geometric methods. The exactness of the BGGL-complex for
arbitrary subsets of simple roots, however, remained open \cite[9.3.19]{b-Kumar02}.

Here we consider two problems related to the BGG-resolution. On the one hand we answer
Kumar's question as anticipated in \cite[3.28]{a-Kumar90}. On the other hand,
we are interested in analogs of both the BGG-resolution and the BGGL-resolution for 
quantized enveloping algebras.
Such analogs appear in various places in the literature \cite{inp-Rosso91}, \cite{a-FF92},
\cite{a-Malikov92}, \cite{phd-Caldero}, \cite[4.5.5]{b-Joseph}.
However, these presentations are either rather short, treat only special cases, or rely on 
specialization arguments. The present note grew out of the authors' desire to obtain
an easily accessible proof of the exactness of a quantum analog of the BGGL-resolution,
avoiding specialization.

The construction of the BGG-complex can be translated into the quantum group setting word by 
word without further ado. The usual proof of the exactness of the BGG-complex, however, 
relies on the standard resolution in Lie algebra homology 
(cp. also \cite[Chapter 9.1]{b-Kumar02}). This standard resolution is written in terms of 
exterior powers of $\gfrak/\pfrak_S$ which are not available in the quantum case 
(note however that the $SL(N)$ case is treated in \cite{inp-Rosso91} in terms of Koszul 
resolutions). 
At this point other authors revert to specialization \cite{a-Malikov92},
\cite{phd-Caldero}. In his book \cite[4.5.5]{b-Joseph} A. Joseph also suggests specialization
or refers to the paper \cite{a-GabJo81} which in a more general setting relies on more 
advanced homological and representation theoretic methods.
 
Here we give an elementary proof of the exactness of the BGG-resolution which does not depend 
on the standard resolution and which also works in the quantum group setting avoiding 
specialization. Hence we prove exactness not only for `almost all $q$' as in 
\cite[Theorem 3.3]{a-Malikov92}, or for all $q$ but certain algebraic numbers
as in \cite[IV.2.3]{phd-Caldero}, but for all $q\in \C\setminus\{0\}$ which are not a root 
of unity. The main idea is to show that the homology group $H_n(C_\ast)$ of the
BGG-complex $(C_\ast,\dif_\ast)$ is an integrable
$U(\gfrak)$-module. Hence for symmetrizable Kac-Moody algebras $H_n(C_\ast)$ vanishes for 
$n\ge 1$. The BGGL-resolution 
$(C_\ast^S,\dif^S_\ast)$ of a parabolic subalgebra $\pfrak_S\subset \gfrak$ is then
constructed as a quotient complex of $(C_\ast,\dif_\ast)$. Its exactness is verified
by showing that the kernel of the canonical surjection $\pi:C_\ast\rightarrow C_\ast^S$
admits a filtration such that the associated graded complex is isomorphic to a direct sum
of BGG-resolutions for the Levi factor $\lfrak_S$ of $\pfrak_S$.

The proof of the exactness of the BGG-complex given here does not work for nonsymmetrizable 
Kac-Moody algebras.
However, in view of Kumar's result \cite{a-Kumar90}, our construction of a quotient complex 
establishes the BGGL-resolution for arbitrary Kac-Moody algebras and arbitrary subsets $S$ 
of the set of simple roots.

\enlargethispage{\baselineskip}
We are very grateful to S.~Kumar for encouraging comments and for pointing out the open
problem \cite[9.3.19]{b-Kumar02} to us. Furthermore we would like to thank 
P.~Caldero and A.~Joseph for providing us with the reference \cite{phd-Caldero}. Thanks are also 
due to H.~Yamane who helped us to understand the equivalence between \cite[3.2.9(vi)]{b-Joseph}
and the quantum Serre relations.

\section{Preliminaries}
In this section we fix the usual notations.
We also give a proof of Lemma \ref{Wlem} which we didn't find in the literature in this form.
It is a returning technical gadget in all later proofs. We mainly refer
to the monograph \cite{b-Kumar02} where original references can be found.

\subsection{Kac-Moody algebras}

Let $\gfrak$ be an arbitrary Kac-Moody algebra over $\C$ with generalized Cartan matrix
$A=(a_{ij})_{1\le i,j\le r}$ satisfying \cite[1.1.1]{b-Kumar02}. Let $(\hfrak,\pi,\pi^\vee)$ be a realization of $A$.
Hence $\hfrak$ is a complex vector space of dimension $r+\mbox{corank}(A)$ and
\begin{align*}
  \pi=\{\alpha_i\}_{1\le i\le r}\subset \hfrak^\ast,\qquad 
  \pi^\vee=\{\alpha_i^\vee\}_{1\le i\le r}\subset \hfrak
\end{align*}
are linearly independent sets satisfying $\alpha_i(\alpha_j^\vee)=a_{ji}$. Recall that 
$\gfrak$ is generated by $\hfrak$ and symbols  $e_i$, $f_i$, $1\le i\le r$, with defining 
relations given for instance in \cite[1.1.2]{b-Kumar02}. Let $\wurz$ denote the set of
roots associated with $(\gfrak,\hfrak, \pi, \pi^\vee)$ and let $\wurz^+\subset\wurz$ be the 
set of positive
roots. For $\alpha\in \wurz$ we write $\gfrak_\alpha$ to denote the corresponding root space,
$\gfrak=\nfrak\oplus\hfrak\oplus \nfrak^-$ for the
triangular decomposition, and $\bfrak=\nfrak\oplus\hfrak$ for the Borel subalgebra.

Let $U(\cdot)$ denote the universal enveloping algebra functor. For $\lambda\in \hfrak^\ast$ 
let $\C_\lambda$ be the one dimensional $\bfrak$-module on which $\hfrak$ acts via $\lambda$
and $\nfrak$ acts trivially, and let $V^\lambda=\ug\ot_{U(\bfrak)}\C_\lambda$
be the corresponding Verma module.

We write $\wlat=\{\lambda\in \hfrak^\ast\,|\, \lambda(\alpha_i^\vee)\in \Z \mbox{ for all }
\alpha_i^\vee \in \pi^\vee\}$ to denote the weight lattice and $\wlat^+\subset\wlat$ to
denote the set of dominant integral weights. For $\lambda\in \wlat^+$
let $V(\lambda)$ be the maximal integrable quotient of $V^\lambda$. 
By \cite[2.1.5, 2.1.7]{b-Kumar02} one has 
\begin{align}\label{Vlamdef}
  V(\lambda)=V^\lambda\big/ \sum_{\alpha_i\in \pi} \U(\gfrak)
  (f_i^{\lambda(\alpha_i^\vee)+1}\ot 1).
\end{align}
The $\U(\gfrak)$-module $V(\lambda)$ is irreducible for symmetrizable $\gfrak$ 
\cite[2.2.6, 3.2.10]{b-Kumar02}.

Let $S\subseteq \pi$ be a set of simple roots. Define $\wurz_S=\wurz\cap\Z S$ and
$\wurz_S^+=\wurz^+\cap\Z S$. The following subspaces of $\gfrak$
are Lie subalgebras:
\begin{align*}
\nfrak_S&=\bigoplus_{\alpha\in \wurz^+_S}\gfrak_\alpha &
\nfrak^-_S& = \bigoplus_{\alpha\in \wurz^+_S}\gfrak_{-\alpha} &
\ufrak_S&=\bigoplus_{\alpha\in \wurz^+\setminus\wurz^+_S}\gfrak_\alpha\\
\lfrak_S&=\nfrak_S\oplus\hfrak\oplus\nfrak^-_{S}&
\pfrak_S&=\lfrak_S\oplus\ufrak_S&
\ufrak_S^-&=\bigoplus_{\alpha\in \wurz^+\setminus\wurz^+_S}\gfrak_{-\alpha}.
\end{align*}
Define $\wlat_S^+=\{\lambda \in \wlat\,|\, \lambda(\alpha_i^\vee)\in \N_0\, 
\mbox{ for all } \alpha_i\in S\}$. 
For $\lambda\in \wlat^+_S$ let $M(\lambda)$ denote the maximal integrable
$\lfrak_S$-module of highest weight $\lambda$. Note that $M(\lambda)$ can be regarded as a 
$\pfrak_S$-module by the given action of $\lfrak_S$ and the trivial action of $\ufrak_S$. 
Let now  $V^{M(\lambda)}=\ug \ot_{U(\pfrak_S)}M(\lambda)$ denote the generalized Verma 
module of highest weight $\lambda$.

\subsection{The Bruhat order}
For  $1\le i \le r$ let $s_i\in \mbox{Aut}(\hfrak^\ast)$ denote the reflection defined by
\begin{align*}
  s_i(\chi)=\chi-\chi(\alpha _i^\vee)\alpha_i \qquad
  \mbox{ for $\chi\in \hfrak^\ast$.}
\end{align*}
The Weyl group $W\subset \mbox{Aut}(\hfrak^\ast)$ is the subgroup generated by 
$\{s_i\,|\,1\le i\le r\}$. Recall that $W$ is a Coxeter group \cite[1.3.21]{b-Kumar02}.
Let $l$ denote the length function on $W$. Let $\rho\in \hfrak^\ast$ be an element such
that $\rho(\alpha_i^\vee)=1$ for all $\alpha_i^\vee\in \pi^\vee$. 
Recall that the shifted action of the Weyl group $W$ on $P$ is defined
in terms of the ordinary Weyl group action by
\begin{align*}
  w.\mu=w(\mu+\rho)-\rho.
\end{align*}
The shifted action does not depend on the choice of $\rho$.
For any subset $S\subset \pi$ let $W_S\subset W$ denote the subgroup generated by the 
reflections corresponding to simple roots in $S$. Moreover, define
\begin{align*}
  W^S=\{w\in W\,|\, l(vw)\ge l(w) \mbox{ for all } v\in W_S\}.
\end{align*}
The following well-known result will be 
frequently used throughout this note.
\begin{lemma}\label{kost-lem}{\upshape \cite[1.3.17]{b-Kumar02}}
  Any element $w\in W$ can be decomposed uniquely in the form $w=w_Sw^S$ where 
  $w_S\in W_S$ and $w^S\in W^S$. Moreover, this decomposition satisfies $l(w)=l(w_S)+l(w^S)$.
\end{lemma}
Let $T=\{ws_iw^{-1}\,|\, \alpha_i\in \pi, w\in W\}$ denote the set of reflections in $W$
corresponding to real roots.
Following \cite{a-Lep77} 
for $w,w'\in W$ we write $w\rightarrow w'$ if there exists
$t\in T$ such that $w=t w'$ and $l(w)=l(w')+1$.
The Bruhat order $\le $ on $W$ is then given by the relation
\begin{align*}
w\le w' \Leftrightarrow \
&\text{there exists } n\ge 1,\ w_2,\ldots ,w_{n-1}\in W,\\
&\text{such that }w=w_1\rightarrow w_2\rightarrow \ldots
\rightarrow w_n=w'.
\end{align*}
\begin{lemma}\label{Wlem}
    Let $w_S,w'_S\in W_S$ and $w,w'\in W^S$ such that relation
    $w_Sw\rightarrow w'_Sw'$ holds.\\
 (i) One has $l(w)\ge l(w')$.\\
 (ii) The relation $l(w)=l(w')$ holds if and only if $w=w'$. In this case one
 additionally has $w_S\rightarrow w'_S$ in $W_S$.
  \end{lemma}
  
  \noindent{\bf Proof:}
  Write $w_S=s_{i_1}\cdots s_{i_k}$ and $w=s_{i_{k+1}}\cdots s_{i_m}$,
  where $k=l(w_S)$, $m-k=l(w)$, and $\alpha _{i_j}\in S$ for $1\le j\le k$.
  Since $l(w_Sw)=l(w_S)+l(w)$ by Lemma~\ref{kost-lem},
  $s_{i_1}\cdots s_{i_m}$ is
  a reduced decomposition of $w_Sw$ into simple reflections.
  By \cite[1.3.16]{b-Kumar02} there exists $j\in
  \{1,2,\ldots ,m\}$ such that
  \begin{align}\label{eq-wdecomp}
    w'_Sw'=s_{i_1}\cdots s_{i_{j-1}}s_{i_{j+1}}\cdots s_{i_m}.
  \end{align}
  Assume first that $j>k$. Then at least the first $k$ factors
  on the right hand side of (\ref{eq-wdecomp}) belong to $W_S$, and hence
  Lemma~\ref{kost-lem} implies that $l(w')=l(w'_Sw')-l(w'_S)\le m-1-k=l(w)-1$.
  On the other hand, if $j\le k$ then $w'_Sw'=
  s_{i_1}\cdots s_{i_{j-1}}s_{i_{j+1}}\cdots s_{i_k}w$ and hence by
  Lemma \ref{kost-lem} one obtains that $w=w'$ and
  $w'_S=s_{i_1}\cdots s_{i_{j-1}}s_{i_{j+1}}\cdots s_{i_k}$.
  $\blacksquare$

\section{The BGGL-resolution as a quotient of the BGG-resolution}

\subsection{The BGG-resolution}
\label{ss-BGGres}

We recall the construction of the Bernstein-Gelfand-Gelfand complex following \cite{a-BGG75}. 
This construction crucially depends on the following result about 
homomorphisms of Verma modules which holds for arbitrary Kac-Moody algebras.
Consult \cite[9.C]{b-Kumar02} for references to the original literature.
\begin{theorem}\label{vermainc}{\upshape \cite[9.2.3]{b-Kumar02}} For any $\mu\in \wlat^+$ one 
  has
  \begin{align*}
    \dim\mbox{\upshape Hom}_{\ug}(V^{w.\mu},V^{w'.\mu})=
      \begin{cases} 1 & \mbox{if $w\le w'$,} \\ 0 &\mbox{else.}
    \end{cases}
  \end{align*}
\end{theorem}
Fix a dominant integral weight $\mu\in \wlat^+$. For all 
$n\in \N_0$ define
\begin{align*}
  C_n:=\bigoplus_{w\in W,\, l(w)=n} V^{w.\mu}.
\end{align*}
One constructs $\ug $-module maps $\dif _n:C_n\rightarrow C_{n-1}$
for all $n\in \N$ as in \cite[$\S$ 10]{a-BGG75}.
More explicitly using the above theorem, for every $w\in W$, fix an embedding 
$V^{w.\mu}\subset V^\mu$. Then, using Theorem \ref{vermainc} again, for all $w,w'\in W$ 
with $w\le w'$ one has a fixed embedding $f_{w,w'}:V^{w.\mu}\rightarrow V^{w'.\mu}$.

A quadruple $(w_1,w_2,w_3,w_4)$ of elements of $W$ is called a square
if $w_2\neq w_3$ and 
\begin{align*}
  w_1\rightarrow w_2\rightarrow w_4, \qquad
  w_1\rightarrow w_3\rightarrow w_4.
\end{align*}
By \cite[Lemma 10.4]{a-BGG75} to each arrow 
$w_1\rightarrow w_2$, where $w_1,w_2\in W$, one can 
assign a number $s(w_1,w_2)=\pm 1$ such that for every square, the
product of the numbers assigned to the four arrows occurring in it is 
$-1$.

The differential $\dif _n:C_n\rightarrow C_{n-1}$ is now defined by
\begin{align*}
  \dif_n=\sum_{\makebox[0cm]{$w,w'\in W,\,l(w)=n\atop w\rightarrow w'$}}
  s(w,w')f_{w,w'}
\end{align*}
where $f_{w,w'}$ is extended by zero to all of $C_n$.
Let moreover $\vep_\mu:V^\mu\rightarrow V(\mu)$ denote the canonical surjection.
By construction the relation $\dif_{n-1}\circ \dif_n=0$ holds. Hence one obtains a complex
\begin{align}\label{eq-BGGres}
   \cdots \longrightarrow C_{2}
    \stackrel{\dif_2}{\longrightarrow} C_1
    \stackrel{\dif_1}{\longrightarrow} C_0
    \stackrel{\vep_\mu}{\longrightarrow} V(\mu)\longrightarrow 0.
  \end{align}

\begin{theorem}\label{BGG-theorem}{\upshape \cite[9.2.20, 9.3.14]{b-Kumar02}}
  The complex (\ref{eq-BGGres}) is exact.
\end{theorem}

In order to generalize this theorem
to Verma modules over quantized enveloping algebras
we are going to give a new proof for the case of symmetrizable Kac-Moody algebras. 
For any $\alpha _i\in \pi $ let $U_i$ denote the subalgebra of
$\ug $ generated by the subspace
$\gfrak _{\alpha _i}\oplus \gfrak _{-\alpha _i}$ of $\gfrak $.
Then $U_i$ is isomorphic to $U(\mathfrak{sl}_2)$ for all $i$.
For any $n\in \N _0$ and $\alpha _i\in \pi $ define subsets $W^+_{n,i}$
and $W^-_{n,i}$ of $W$ by
\begin{align}
\begin{aligned}
W^+_{n,i}&:=\{w\in W\,|\,l(w)=n,\,l(s_iw)=n+1\},\\
W^-_{n,i}&:=\{w\in W\,|\,l(w)=n,\,l(s_iw)=n-1\}.
\end{aligned}
\end{align}
Note that $w\in W^+_{n,i}$ if and only if $w^{-1}\alpha_i>0$, or equivalently if and only if
$w^{-1}\alpha_i^\vee\in \sum_{i=1}^r\N_0\alpha_i^\vee$.
The following well-known lemma which follows from \cite[Lemma~4.2.7]{b-Joseph} will be 
needed in the proof of Theorem
\ref{BGG-theorem} for symmetrizable Kac-Moody algebras.

\begin{lemma}\label{l-integrability}
Let $\gfrak$ be an arbitrary Kac-Moody algebra.
For all $w\in W^+_{n,i}$ and $\mu \in P^+$ the $\ug $-module $V^{w.\mu }/
V^{s_iw.\mu }$ is locally finite under the action of $U_i$.
\end{lemma}

The exactness of the sequence (\ref{eq-BGGres}) at $V(\mu)$ is just the definition of
$V(\mu )$ as a quotient of the Verma module $V^\mu$.
The exactness at $C_0$ also holds by the definition (\ref{Vlamdef}) of $V(\mu)$.
For arbitrary Kac-Moody Lie algebras it is not known if the module $V(\mu)$ is simple.
In the symmetrizable case, however, it suffices to show that $V(\mu)$ is integrable in order
to obtain that $V(\mu)$ is simple \cite[Proposition 4.2.8(iii)]{b-Joseph}. 
We use the same idea to verify exactness at $C_n$ for $n\ge 1$.

\begin{proposition}\label{integrable} Let $\gfrak$ be an arbitrary Kac-Moody algebra.
  For any $n\in \N$ the homology group $H_n(C_\ast)$ of the 
  complex (\ref{eq-BGGres}) is an integrable $\ug$-module.
\end{proposition}

\noindent{\bf Proof:}
Take $\alpha _i\in \pi $. We have to show that $\ke (\dif _n)/\im (\dif _{n+1})$ is locally 
finite under the action of $U_i$.
It follows from Lemma \ref{Wlem} for $S=\{\alpha _i\}$ that for all
$w\in W^+_{n,i}$ one has $\dif _{n+1}(V^{s_iw.\mu })\subset V^{w.\mu }
\oplus \bigoplus _{w'\in W^-_{n,i}}V^{w'.\mu }$. 
Hence one obtains an isomorphism of $\ug$-modules
\begin{align*}
  C_n\Bigg/\Bigg(\bigoplus _{w'\in W^-_{n,i}}V^{w'.\mu }+\sum _{w\in W^+_{n,i}}
\dif _{n+1}(V^{s_iw.\mu })\Bigg)\cong \bigoplus_{w\in W^+_{n,i}}V^{w.\mu}/V^{s_iw.\mu}.
\end{align*}
By Lemma \ref{l-integrability} the right hand side is locally finite with respect to 
$U_i$. Hence, so is
\begin{align}\label{express}
  \ke(\dif_n)\Bigg/\Bigg(\Bigg(\ke(\dif_n)\cap\bigoplus _{w'\in W^-_{n,i}}V^{w'.\mu }\Bigg)
    +\im(\dif_{n+1})\Bigg).
\end{align}
Again, for $S=\{\alpha_i\}$ in Lemma \ref{Wlem}, one obtains for any 
$w',w''\in W^-_{n,i}$ that
$w'\rightarrow s_iw''$ if and only if $w'=w''$. Hence the restriction
of $\dif _n$ to the subspace $\bigoplus _{w'\in W^-_{n,i}}V^{w'.\mu }$
of $C_n$ is injective and the expression (\ref{express}) coincides with $H_n(C_\ast)$.
$\blacksquare$

\vspace{.5cm}

\noindent
{\bf Proof of Theorem \ref{BGG-theorem}} (for symmetrizable Kac-Moody algebras){\bf:}
By the above proposition the homology $H_n(C_\ast)$ is integrable and therefore
isomorphic to a direct sum of simple $\ug$-modules \cite[2.2.7]{b-Kumar02}.
These simple modules are subquotients of $V^\mu $ and hence isomorphic to
$V(\mu )$ \cite[2.1.16, 2.2.4]{b-Kumar02}. However $C_n$ does not
contain vectors of weight $\mu $ for $n\ge 1$. Thus $H_n(C_\ast)=0$ holds
also for $n\ge 1$.
$\blacksquare$

\subsection{A quotient complex}
\label{ss-quotcomplex}

Define
\begin{align*}
  C_n^S=\bigoplus_{w\in W^S,\, l(w)=n} V^{M(w.\mu)}
\end{align*}
and let $\pi_n:C_n\rightarrow C_n^S$ denote the canonical $U(\gfrak)$-module projection.
More explicitly, $\pi_n=\bigoplus_{w\in W,\,l(w)=n} \pi_{n,w}$ where
$\pi_{n,w}:V^{w.\mu}\rightarrow V^{M(w.\mu)}$ for $w\in W^S$ denotes the canonical
$U(\gfrak)$-module projection and $\pi_{n,w}:V^{w.\mu}\rightarrow 0$ if $w\notin W^S$.
Note that by definition of $M(w.\mu)$ one has 
$\ke (\pi_{n,w})=\sum_{\alpha_i\in S} f_{s_iw,w}(V^{s_iw.\mu})$ if $w\in W^S$ and 
hence
\begin{align}\label{kerpin}
  \ke(\pi_n) =\bigoplus_{w\notin W^S,\, l(w)=n}V^{w.\mu}\bigoplus_{w\in W^S,\,l(w)=n}
  \sum_{\alpha_i\in S} f_{s_iw,w}(V^{s_iw.\mu}).
\end{align}
\begin{lemma}
  For all $n\in \N$ one has $\dif_n(\ke (\pi_n))\subseteq \ke (\pi_{n-1})$.
  Moreover, $\vep_\mu(\ke(\pi_0))=0$.
\end{lemma}

\noindent{\bf Proof:} By \cite[9.2.14]{b-Kumar02} for all $w,w'\in W^S$, $w\rightarrow w'$,
the inclusion $f_{w,w'}:V^{w.\mu}\rightarrow V^{w'.\mu}$ maps
$\sum_{\alpha_i\in S} f_{s_iw,w}(V^{s_iw.\mu})$ to 
$\sum_{\alpha_i\in S} f_{s_iw',w'}(V^{s_iw'.\mu})$.

Consider now $\dif_n$ restricted to $V^{w.\mu}$ for $w\notin W^S$. 
By Lemma \ref{kost-lem} one can write $w=w_Sw^S$ with 
$w_S\in W_S\setminus\{e\}$, $w^S\in W^S$,
and $l(w)=l(w_S)+l(w^S)$. If $l(w_S)\ge 2$ then
$\dif_n(V^{w.\mu})\subseteq \bigoplus_{w'\notin W^S}V^{w'.\mu}$ by Lemma \ref{Wlem}.
If $l(w_S)=1$ then $w_S=s_i$ for some $\alpha_i\in S$ and hence by Lemma 
\ref{Wlem} one obtains
\begin{align*}
  \dif_n(V^{w.\mu})\subseteq \bigoplus_{w'\notin W^S}V^{w'.\mu}
  \bigoplus f_{s_iw^S,w^S}(V^{s_iw^S.\mu})\subseteq \ke (\pi_{n-1})
\end{align*}
also in this case.
$\blacksquare$

\vspace{.5cm}

\noindent By the above lemma one obtains a short exact sequence of complexes
\begin{align*}
  0\rightarrow \ke (\pi_\ast)\rightarrow C_\ast\stackrel{\pi}{\rightarrow} C_\ast ^S
  \rightarrow 0.
\end{align*}
Let $\dif ^S_n:C_n^S\rightarrow C_{n-1}^S$ denote the induced differentials on $C^S_\ast$.
Note that by construction $(C_\ast ^S,\dif^S_\ast)$ coincides with the complex
considered in \cite[9.2.17]{b-Kumar02}. 

\subsection{Exactness}
We now give a generalization of the BGG-resolution for arbitrary Kac-Moody algebras and
arbitrary subsets $S\subset \pi$. The proof uses Theorem \ref{BGG-theorem}.
In the nonsymmetrizable case one only has Kumar's original proof \cite{a-Kumar90} of
Theorem \ref{BGG-theorem} which relies on more involved geometric methods. 
In this case the following result was anticipated in \cite[3.28]{a-Kumar90}.

\begin{theorem}\label{Lepowsky-theorem}
  The complex $(C_\ast ^S,\dif^S_\ast)$ is exact.
\end{theorem}

\noindent{\bf Proof:}
Using the long exact homology sequence it suffices to show that
the complex $(\ke (\pi_\ast),\dif_\ast|_{\ke (\pi_\ast)})$ is exact.

Consider the $\N _0$-filtration $\cF $ of $C_*$ defined by
\begin{align}\label{eq-filtration}
\cF _k(C_n)=
\bigoplus _{w_S\in W_S,w^S\in W^S\atop l(w_Sw^S)=n,\,l(w^S)\le k}V^{w_Sw^S.\mu }.
\end{align}
Lemma \ref{Wlem} implies that $(C_*,\dif _*)$ is a filtered complex.
For $w^S\in W^S$ with $l(w^S)=k$ define
\begin{align*}
D_{n,w^S}=\begin{cases}
\displaystyle
\bigoplus _{w_S\in W_S,\,l(w_S)=n-k}V^{w_Sw^S.\mu } & \text{if $k\not=n$,}\\
\sum _{\alpha _i\in S}f_{s_iw^S,w^S}(V^{s_iw^S.\mu }) & \text{if $k=n$.}
\end{cases}
\end{align*}
By (\ref{kerpin}) the associated graded complex of $(\ke (\pi _*),\dif |_{\ke (\pi _*)})$
has homogeneous components
$\gr ^\cF _k\ke (\pi _n)\cong \bigoplus _{w^S\in W^S,\,l(w^S)=k}D_{n,w^S}$.
Moreover, by Lemma \ref{Wlem} one obtains that $D_{*,w^S}$ is a subcomplex of the graded 
complex $\gr^\cF \ke (\pi _*)$ for every $w^S\in W^S$.
We complete the proof of the theorem by showing that $D_{*,w^S}$ is exact.

To this end recall that for any $\lambda\in P_S^+$
one has an isomorphism of $\hfrak $-modules
\begin{align}\label{eq-Vzerl}
  V^\lambda\cong \uusm\ot \unsm v_\lambda
\end{align}
where $\hfrak$ acts on the right hand side diagonally and $v_\lambda\in V^\lambda$ denotes 
a highest weight vector.
With respect to this factorization the inclusion 
$f_{w_Sw^S,w_S'w^S}:V^{w_Sw^S.\mu}\hookrightarrow V^{w_S'w^S.\mu}$ satisfies
\begin{align}\label{f=idotg}
  f_{w_Sw^S,w_S'w^S}=\id\ot g_{w_S,w_S'}
\end{align}
where $g_{w_S,w_S'}:\unsm v_{w_Sw^S.\mu}\rightarrow \unsm v_{w_S'w^S.\mu}$
denotes the inclusion of $U(\lfrak_S)$-modules.
Note moreover that for  
\begin{align*}
  D^S_{j,w^S}:=\bigoplus_{w_S\in W_S,\,l(w_S)=j}\unsm v_{w_Sw^S.\mu}
\end{align*}  
and 
\begin{align*}
  \delta^S_{j,w^S}:=\sum_{\makebox[1cm]{$w_S, w_S'\in W_S,\atop
 l(w_S)=j,\,w_S\rightarrow w_S'$}}
                    s(w_Sw^S,w_S'w^S) g_{w_S,w_S'}
\end{align*}  
the pair $(D^S_{\ast ,w^S},\delta^S_{\ast ,w^S})$ is a complex isomorphic
to the BGG-resolution
of the $U(\lfrak_S)$-module $M(w^S.\mu)$
and therefore exact by Theorem \ref{BGG-theorem}.
By (\ref{eq-Vzerl}) and (\ref{f=idotg}) up to a shift of degree the complex
$D_{*,w^S}$ is isomorphic to $\uusm \otimes \bar{D}^S_{*,w^S}$,
where $\bar{D}^S_{*,w^S}$ is the right truncation
$$
  \cdots \rightarrow D^S_{3,w^S} \rightarrow D^S_{2,w^S}
\rightarrow D^S_{1,w^S}\rightarrow \delta ^S_{1,w^S}(D^S_{1,w^S})\rightarrow 0
$$
of the resolution $D^S_{*,w^S}$, and hence it is exact.
$\blacksquare$

\subsection{The quantum BGGL-resolution}

In this section we discuss the quantum versions of Proposition \ref{integrable} and Theorems
\ref{BGG-theorem} and \ref{Lepowsky-theorem}.
Let $\uqg $ denote the quantized enveloping algebra of a symmetrizable Kac-Moody
algebra $\mathfrak{g}$. The main advantage of our approach in this paper is
that the following remarks are valid for both, the definition of $\uqg $ in
\cite[3.2.9]{b-Joseph} over the field $k(q)$,
where $k$ is a field of characteristic zero,
and the definition in \cite[6.1.2]{b-KS} over the field of complex numbers,
where $q\in \mathbb{C}\setminus \{0\}$ is not a root of unity. Note that in \cite{b-KS} the
algebra $\uqg$ is defined only for semisimple $\gfrak$. However for a symmetrizable Kac-Moody
algebra one can consider the same presentation in terms of generators and relations.
Note that the equivalence of the relations \cite[3.2.9(vi)]{b-Joseph} with the quantum
Serre relations follows from \cite[4.1.17]{b-Joseph} analogously to the proof of 
\cite[Theorem 9.11]{b-Kac1}.

Along the lines of Section \ref{ss-BGGres} one obtains a sequence $C_\ast$ of
$\uqg $-modules analogous to (\ref{eq-BGGres}). To this end one has to
verify that the relation
\begin{align*}
\dim \Hom _{\uqg }(V^{w.\mu },V^{w'.\mu })= 1\qquad \mbox{ if } w\le w'
\end{align*}
holds for all $\mu\in P^+$. This can be seen following 
\cite[4.4.7, 4.4.15]{b-Joseph}. The proof also holds for $q\in \mathbb{C}\setminus \{0\}$ 
not a root of unity.
Now one can take the proof of Proposition \ref{integrable} literally to
obtain the integrability of the homology $H_n(C_*)$. The necessary
quantum analog of Lemma \ref{l-integrability} can be found in
\cite[4.3.5]{b-Joseph}. The proof given there works for $q\in \mathbb{C}\setminus \{0\}$ 
not a root of unity.

Note that by \cite[4.3.6(i)]{b-Joseph} which implies the second half of
\cite[Theorem 4.3.10]{b-Joseph} and by \cite[3.4.9]{b-Joseph} the proof of Theorem 
\ref{BGG-theorem} also holds for $\uqg$. Again the argument also works for
$q\in \mathbb{C}\setminus \{0\}$ not a root of unity.

Note that \cite[9.2.14]{b-Kumar02} translates into the quantum setting. Hence
one may construct quantum analogs of the complexes $\ke (\pi _*)$
and $C^S_*$ as in Section \ref{ss-quotcomplex}. Finally, the proof of Theorem 
\ref{Lepowsky-theorem} also applies literally in the quantum setting.

%\bibliographystyle{amsalpha}
%\bibliography{litbank2}

\providecommand{\bysame}{\leavevmode\hbox to3em{\hrulefill}\thinspace}
\providecommand{\MR}{\relax\ifhmode\unskip\space\fi MR }
\providecommand{\MRhref}[2]{%
  \href{http://www.ams.org/mathscinet-getitem?mr=#1}{#2}
}
\providecommand{\href}[2]{#2}

\pagebreak

\textsc{Istv\'an Heckenberger, Mathematisches Institut, Universit\"at Leipzig,
         Augustusplatz 10-11, 04109 Leipzig, Germany.}
       
 \textit{E-mail address:} \texttt{heckenberger@math.uni-leipzig.de}

\vspace{.5\baselineskip}
                                                
\textsc{Stefan Kolb, Mathematics Department, Virginia Polytechnic Institute
         and State University, Blacksburg, VA 24061, USA.}
         
 \textit{E-mail address:} \texttt{kolb@math.vt.edu}

\end{document}